\input amstex
\documentstyle {amsppt}
\NoBlackBoxes
\leftheadtext{Igor Belegradek}
\rightheadtext{Counting negatively curved manifolds}
\topmatter
\title
{Counting open negatively curved manifolds up to tangential
homotopy equivalence}
\endtitle
\author Igor Belegradek 
\endauthor
\address
Department of Mathematics,
University of Maryland, College Park, MD 20742
\endaddress
\email
{\rm igorb\@math.umd.edu}
\endemail
\date
May 14, 1998
\enddate
\abstract
Under mild assumptions on a group $\pi$,
we prove that the class  
of complete Riemannian $n$--manifolds
of uniformly bounded negative sectional curvatures and 
with the fundamental groups isomorphic to $\pi$
breaks into finitely many tangential homotopy types.
It follows that many aspherical manifolds do not admit
complete negatively curved metrics with
prescribed curvature bounds. 
\endabstract 
\endtopmatter
\document

\heading
\S 1.~Introduction
\endheading

This paper is an attempt to understand the topology
of complete infinite volume Riemannian manifolds of 
negative sectional curvature.
Every smooth open manifold admits a (possibly incomplete)
Riemannian metric of
negative sectional curvature [23].
However, the universal cover of a complete negatively curved 
manifold is diffeomorphic to the Euclidean space,
hence the homotopy type of a complete negatively curved manifold
is determined by its fundamental group.
 
For closed 
(or, more generally, finite volume complete)
negatively curved manifolds, the fundamental group seems to
encode most of the topological information.
By contrast, complete
negatively curved manifolds of {\it infinite} volume
with isomorphic fundamental groups
may be very different topologically.
For example, the total space of any vector bundle
over a closed negatively curved manifold admits
a complete Riemannian metric of sectional curvature 
pinched between two negative constants [1].

Given a group $\pi$, 
a positive integer $n$, and real numbers $a\le b<0$,
consider the class $\Cal M_{a,b,\pi, n}$
of $n$--manifolds with fundamental groups
isomorphic to $\pi$ that can be given complete Riemannian 
metrics of sectional
curvatures within $[a,b]$.

In this paper we discuss 
under what assumptions on $\pi$ the class
$\Cal M_{a, b, \pi, n}$ breaks into {\it finitely
many} tangential homotopy types.
Recall that a homotopy equivalence of 
smooth manifolds of the same dimension $f: N\to L$
is called tangential if the vector bundles
$f^*TL$ and $TN$ are stably isomorphic.
For example, any map that 
is homotopic to a diffeomorphism is a tangential 
homotopy equivalence. 
Note that, unless $N$ is a closed manifold, 
the bundles $f^*TL$ and $TN$ are isomorphic
iff they are stably isomorphic.

If the cohomological dimension of 
$\pi$ is equal to $n$, 
every manifold from the class 
$\Cal M_{a,b,\pi, n}$ is closed.
In that case, in fact, $\Cal M_{a,b,\pi, n}$ 
falls into finitely many diffeomorphism classes [5].

Here is a way to produce infinitely many homotopy equivalent
negatively curved manifolds of the same dimension
such that no two of them are tangentially
homotopy equivalent. 
Let $M$ be a closed negatively curved manifold
with $H_{4m}(M,\Bbb Q)\neq 0$ for some $m>0$.
(Such examples abound. For instance, any closed
complex hyperbolic or quaternion hyperbolic manifold 
is such.
Most arithmetic closed real hyperbolic orientable 
manifolds have nonzero Betti numbers in all dimensions [30].)
Then, by an elementary K-theoretic argument,
for any $k\ge\dim (M)$ there exists a
infinite sequence of rank $k$ vector bundles over $M$
such that no two of them have 
tangentially homotopy equivalent total spaces. 
Yet, thanks to a theorem of Anderson [1], 
the total space of any vector bundle over $M$
can be given a negatively curved metric.

The following is a simplified form of our main result.

\proclaim{Theorem~1.1}
Let $\pi$ be the fundamental group of a finite aspherical complex. 
Suppose that $\pi$ is not virtually nilpotent
and that $\pi$ does not split as a 
nontrivial amalgamated product
or an HNN-extension over a virtually nilpotent group.

Then, for any positive integer $n$ and any negative reals
$a\le b$, the class $\Cal M_{a, b, \pi, n}$ 
breaks into finitely
many tangential homotopy types.
\endproclaim

Furthermore, by an easy argument 
we can generalize the theorem~$1.1$ to certain
amalgamated products and HNN-extensions. 
For example, suppose that $\pi$ is the fundamental group of
a finite graph of groups such that the edge groups
are virtually nilpotent groups of cohomological
dimension $\le 2$. 
(Note that if $\Cal M_{a, b, \pi, n}\neq\emptyset$,
then any nontrivial virtually nilpotent subgroup of $\pi$ 
of cohomological dimension $\le 2$ is isomorphic to $\Bbb Z$, 
$\Bbb Z\times\Bbb Z$, or the fundamental group
of the Klein bottle.)
Fix $n$ and $a\le b<0$ and suppose that, 
for each vertex group $\pi_v$,
the class $\Cal M_{a, b, \pi_v, n}$ breaks into finitely
many tangential homotopy types. 
Then so does the class $\Cal M_{a, b,\pi, n}$. 

Applying a powerful accessibility result of Delzant
and Potyagailo [16], we deduce the following.

\proclaim{Corollary~1.2}
Let $\pi$ be the fundamental group of a finite aspherical complex.
Assume that any 
nilpotent subgroup of $\pi$ has cohomological dimension $\le 2$.

Then, for any positive integer $n$ and any negative reals
$a\le b$, the class $\Cal M_{a, b, \pi, n}$ 
breaks into finitely
many tangential homotopy types.
\endproclaim

Any torsion free word-hyperbolic group is the fundamental group of
a finite aspherical cell complex [15, 5.24]. 
Moreover, any virtually nilpotent subgroup
of $\pi$ is either trivial or infinite cyclic.

\proclaim{Corollary~1.3}
Let $\pi$ be a word-hyperbolic group.
Then, for any positive integer $n$ and any negative reals
$a\le b$, the class $\Cal M_{a, b, \pi, n}$ 
breaks into finitely
many tangential homotopy types.
\endproclaim

It is worth mentioning that 
in the locally symmetric case a different (and elementary)
argument yields the following.

\proclaim{Theorem~1.4}
Let $\pi$ be a finitely presented torsion--free group
and let $X$ be a nonpositively curved symmetric space.

Then the class of manifolds of the form $X/\rho(\pi)$,
where $\rho\in\text{Hom}(\pi, \text{Isom}(X))$ is a faithful 
discrete representation, falls into finitely many
tangential homotopy types.
\endproclaim

Gromov posed a question (see [1]) whether an (open) thickening
of a finite aspherical cell complex with word-hyperbolic
fundamental group admits a complete negatively curved metric.
The following corollary shows that most thickenings
do {\it not} carry complete negatively curved metric
with {\it prescribed} curvature bounds.

\proclaim{Corollary~1.5} Let $K$ be a finite, connected, 
aspherical cell complex. Let $a\le b<0$ and 
$n>\max\{4, 2\dim (K)\}$.
Suppose that the class $\Cal M_{a, b, \pi_1(K), n}$ 
breaks into finitely
many tangential homotopy types.  

Then 
the set of diffeomorphism classes of open thickenings of $K$
that belong to the class $\Cal M_{a, b, \pi_1(K), n}$ is \bf finite.
\endproclaim

The corollary follows from the fact that
any two tangentially homotopy equivalent thickenings 
of sufficiently high dimension are diffeomorphic 
[29, pp.226--228].
In general, the set of diffeomorphism classes of open $n$-dimensional
thickenings of $K$ is infinite provided $n>\max\{4, 2\dim (K)\}$
and $\oplus_k H_{4k}(K,\Bbb Q)\neq 0$.
Thus, for such $n$ and $K$,
most open $n$-dimensional thickenings of $K$
do not belong to the class $\Cal M_{a, b, \pi_1(K), n}$.
Combining $1.3$ and $1.5$, we deduce the following.

\proclaim{Corollary~1.6} Let $M$ be a smooth aspherical manifold 
such that $\pi_1(M)$ is word-hyperbolic. Let $a\le b<0$ and 
$n>2\dim (M)$.  

Then the set of isomorphism classes of vector bundles over $M$
whose total spaces belong to $\Cal M_{a, b, \pi_1(K), n}$ is \bf finite.
\endproclaim

Actually, we show in [6] that, in case $\dim(M)\ge 3$,
the assumption $n>2\dim (M)$
of the corollary~$1.6$ is redundant. 

\subheading{Synopsis of the paper}
In the section $2$ we prove main lemmas from the 
global Riemannian geometry.
The $3$rd section provides some K-theoretic background.
In the $4$th section we define some invariant of actions.
The $5$th section is devoted to our main results
including the theorem~$1.1$.
In the section $6$ we employ some group theory to deduce
the corollary~$1.2$.
The section $7$ deals with applications to thickenings. 
Applications  to convex-cocompact groups are deduced 
in the section $8$. 
The theorem $1.4$ is proved in the $9$th section. 

\subheading{Acknowledgments}
I am grateful to Werner Ballmann, Mladen Bestvina, 
Thomas Delzant, Martin J.~Dunwoody, Lowell E.~Jones,
Misha Kapovich, Vitali Kapovitch, Bernhard Leeb, John J.~Millson,
Igor Mineyev, James A.~Schafer, Jonathan M.~Rosenberg and
Shmuel Weinberger for helpful discussions and communications.

Special thanks are due to my advisor 
Bill Goldman for his constant interest and support.

\head
\S 2.~Two types of convergence.
\endhead

By an {\it action} of an abstract group $\pi$ on a space $X$
we mean a group homomorphism $\rho:\pi\to\text{Homeo}(X)$.
An action $\rho$ is {\it free} if $\rho(\gamma)(x)\neq x$
for all $x\in X$ and all $\gamma\in \pi\setminus\{\text{id}\}$.
In particular, if $\rho$ is a free action, then $\rho$ is injective. 

\subheading{2.1.~Equivariant pointed Lipschitz topology}
Let $\Gamma_k$ be a discrete
subgroup of the isometry group
of a complete Riemannian manifold $X_k$
and $p_k$ be a point of $X_k$.
The class of all such triples $\{(X_k, p_k, \Gamma_k)\}$ 
can be given the so-called 
equivariant pointed Lipschitz topology [20]; 
when $\Gamma_k$ is
trivial this reduces to the usual pointed Lipschitz topology.
For convenience of the reader
we give here some definitions borrowed from [20].

For a group $\Gamma$ acting on a pointed metric space $(X,p,d)$
the set $\{\gamma\in\Gamma: d(p, \gamma(p))<r\}$ is denoted by $\Gamma(r)$.
An open ball in $X$ of radius $r$ with center at $p$ is denoted by
$B_r(p,X)$.

For $i=1,2$, let $(X_i, p_i)$ be a pointed
complete metric space with the distance function
$d_i$ and let $\Gamma_i$ be a discrete group of isometries
of $X_i$. In addition, assume that $X_i$ is a $C^\infty$--manifold.
Take any $\epsilon>0$.

Then a quadruple $(f_1, f_2, \phi_1, \phi_2)$ of maps 
$f_i:B_{1/\epsilon}(p_i, X_i)\to B_{1/\epsilon}(p_{3-i}, X_{3-i})$
and $\phi_i:\Gamma_i(1/3\epsilon)\to\Gamma_{3-i}$
is called an $\epsilon$--{\it Lipschitz approximation}
between the triples
$(X_1,p_1,\Gamma_1)$ and $(X_2,p_2,\Gamma_2)$
if the following seven condition hold:

$\bullet$ $f_i$ is a diffeomorphism onto its image;

$\bullet$ for each $x_i\in B_{1/3\epsilon}(p_i, X_i)$
and every $\gamma_i\in\Gamma_i(1/3\epsilon)$,  
$f_i(\gamma_i(x_i))=\phi_i(\gamma_i)(f_i(x_i))$;

$\bullet$ for every $x_i, x_i^\prime\in B_{1/\epsilon}(p_i, X_i)$,
$e^{-\epsilon}<d_{3-i}(f_i(x_i),f_i(x_i^\prime))/d_i(x_i,x_i^\prime)<e^\epsilon$; 

$\bullet$ $f_i(B_{1/\epsilon}(p_i, X_i))\supset
B_{(1/\epsilon)-\epsilon}(p_{3-i}, X_{3-i})$ and 
$\phi_i(\Gamma_i(1/3\epsilon))\supset\Gamma_{3-i}(1/3\epsilon-\epsilon)$;

$\bullet$ $f_i(B_{(1/\epsilon)-\epsilon}(p_i, X_i))\supset
B_{1/\epsilon}(p_{3-i}, X_{3-i})$ and 
$\phi_i(\Gamma_i(1/3\epsilon-\epsilon))\supset\Gamma_{3-i}(1/3\epsilon)$;

$\bullet$ $f_{3-i}\circ f_i|_{B_{(1/\epsilon)-\epsilon}(p_i, X_i)}=\text{id}$
and $\phi_{3-i}\circ\phi_i|_{\Gamma_i(1/3\epsilon-\epsilon)}=\text{id}$;

$\bullet$ $d_{3-i}(f_i(p_i),p_{3-i})<\epsilon$.

We say a sequence of triples $(X_k, p_k, \Gamma_k)$
{\it converges} to $(X, p, \Gamma)$ in the 
equivariant pointed Lipschitz topology
if for any $\epsilon>0$
there is $k(\epsilon)$ such that for all $k>k(\epsilon)$,
there exists an $\epsilon$--Lipschitz approximation
between $(X_k, p_k, \Gamma_k)$ and $(X, p, \Gamma)$.

Notice that if all the groups $\Gamma_k$ are trivial,
then $\Gamma$ is trivial; in this case 
we say that that $(X_k,p_k)$ converges to $(X,p)$ in the pointed
Lipschitz topology. 
Note that if $X_k$ is a complete
Riemannian manifold for all $k$,
then the space $X$ is necessarily a $C^\infty$--manifold with 
a complete $C^{1,\alpha}$--Riemannian metric [22]. 

\subheading{Remark~2.2}
For those with the Kleinian groups background
we note that the equivariant pointed Lipschitz topology 
is closely related to the so-called
Chabauty topology [7][14].
Since the Kleinian group theory is an important source of 
examples, we describe the precise relation.
Let $X$ be a complete Riemannian manifold
(e.g.~a hyperbolic space)
with the isometry group $G$.
The set of discrete subgroups of $G$
has the Chabauty topology induced by
a natural topology on the set of closed subsets of $G$, 
namely a sequence of closed subsets $S_i$ converges to $S$
if for any compact $K\subset G$,
the sets $S_i\cap K$ converge to $S\cap K$
in the Hausdorff topology on $K$.
Consider the product
topology on the set of pairs $(\Gamma, p)$
where $\Gamma$ is a discrete subgroup of $G$
and $p\in X$.  
This product topology is in fact equivalent
to the equivariant pointed Lipschitz topology
where $(\Gamma, p)$ is thought of as a triple
$(X,\Gamma, p)$.   

\subheading{2.3.~Pointwise convergence topology}
Suppose that, for some $p_k\in X_k$, 
the sequence  $(X_k, p_k)$ converges to $(X,p)$
in the pointed Lipschitz topology
i.e., for any $\epsilon>0$
there is $k(\epsilon)$ such that for all $k>k(\epsilon)$
there exists an $\epsilon$--Lipschitz approximation $(f_k, g_k)$
between $(X_k, p_k)$ and $(X, p)$.
(Note that if each $X_k$ is a Hadamard manifold,
then for any $p_k\in X_k$, the sequence $(X_k,p_k)$
is precompact in the pointed Lipschitz topology
because the injectivity radius of $X_k$ at $p_k$
is uniformly bounded away from zero [20].)
We say that a sequence $x_k\in X_k$ {\it converges} to
$x\in X$ if for some $\epsilon$
$$d(f_k(x_k), x)\to 0\ \ \text{as}\ \ k\to\infty$$
where $d(\cdot,\cdot)$ is the distance function on $X$
and $f_k$ comes from the $\epsilon$--Lipschitz approximation
$(f_k, g_k)$ between $(X_k, p_k)$ and $(X, p)$.
Trivial examples: if $(X_k, p_k)$ converges to $(X,p)$
in the pointed Lipschitz topology, then $p_k$ converges to $p$;
furthermore, if $x\in X$, the sequence $g_k(x)$ converges to $x$.

Given a sequence of isometries $\gamma_k\in\text{Isom}(X_k)$
we say that $\gamma_k$ {\it converges}, if for any 
sequence $x_k\in X_k$ that converges to $x\in X$, 
$\gamma_k(x_k)$ converges. 
The limiting transformation $\gamma$
that takes $x$ to the limit of $\gamma_k(x_k)$
of $X$ is necessarily an isometry.
Furthermore, if $\gamma_k$ and $\gamma_k^\prime$
converge to $\gamma$ and $\gamma^\prime$ respectively,
then $\gamma_k\cdot\gamma_k^\prime$ converges
to $\gamma\cdot\gamma^\prime$.
In particular, $\gamma_k^{-1}$ converges to $\gamma^{-1}$
since the identity maps $\text{id}_k:X_k\to X_k$
converge to $\text{id}:X\to X$.

Let $\rho_k:\pi\to\text{Isom}(X_k)$ 
be a sequence of isometric actions of a group $\pi$ on $X_k$.
We say that a sequence of actions $(X_k, p_k, \rho_k)$ 
{\it converges in the pointwise convergence topology} 
to an action $(X, p, \rho)$ if 
$\rho_k(\gamma)$ converges to $\rho(\gamma)$
for every $\gamma\in\pi$.
The limiting map $\rho:\Gamma\to\text{Isom}(X)$ that takes 
$\gamma$ to the limit of $\rho_k(\gamma)$ is necessarily
a homomorphism. 
If $\pi$ is generated by a finite set $S$, then
in order to prove that $\rho_k$ converges in the pointwise
convergence topology it suffices to check that 
$\rho_k(\gamma)$ converges, for every $\gamma\in S$. 

A sequence of actions $(X_k, p_k, \rho_k)$
is called {\it precompact in the pointwise convergence
topology} if every subsequence of
$(X_k, p_k, \rho_k)$ has a subsequence that
converges in the pointwise
convergence topology.

Repeating the proof in [28, 4.7], 
it is easy to check that 
a sequence of isometries $\gamma_k\in\text{Isom}(X_k)$ 
has a converging subsequence if, 
for some converging sequence $x_k\in X_k$, 
the sequence $d_k(x_k,\gamma_k(x_k))$
is bounded 
(where $d_k(\cdot,\cdot)$ is the distance function on $X_k$).

Suppose that $\pi$ is a countable group and assume that
for each $\gamma\in\pi$
the sequence $d_k(p_k,\rho_k(\gamma)(p_k))$
is bounded. Then $(X_k, p_k, \rho_k)$
is precompact in the pointwise convergence topology. 
(Indeed, let $\gamma_1\dots\gamma_n\dots$ be the list of 
all elements of $\pi$. 
Take any subsequence $\rho_{k,0}$ of $\rho_k$. 
Pass to subsequence $\rho_{k,1}$ of $\rho_{k,0}$ so that
$\rho_{k,1}(\gamma_1)$ converges.
Then pass to subsequence $\rho_{k,2}$ of $\rho_{k,1}$
such that $\rho_{k,2}(\gamma_2)$ converges, etc.
Then $\rho_{k,k}(\gamma_n)$ converges for every $n$.)

Note that if $\pi$ is generated by a finite set $S$, then
to prove that $\rho_k$ is precompact
it suffices to check that 
$d_k(p_k,\rho_k(\gamma)(p_k))$ is bounded, 
for all $\gamma\in S$ because it implies that 
$d_k(p_k,\rho_k(\gamma)(p_k))$ is bounded, 
for each $\gamma\in\pi$.

\subheading{2.4.~Motivating example}
Let $X$ be a complete Riemannian manifold.
Consider the isometry group
$\text{Isom}(X)$ of $X$ and let $\pi$ be a group. 
The space $\text{Hom}(\pi,\text{Isom}(X))$ has a natural
topology (which is usually called ``algebraic topology''
or ``pointwise convergence topology''), namely
$\rho_k$ is said to converge to $\rho$ if, for each $\gamma\in\pi$,
$\rho_k(\gamma)$ converges to $\rho(\gamma)$ in the Lie group 
$\text{Isom}(X)$.  
Note that if $\pi$ is finitely generated, this topology
on $\text{Hom}(\pi,\text{Isom}(X))$
coincide with the compact-open topology.
Certainly, for any $p\in X$, 
the constant sequence $(X,p)$ converges to itself
in pointed Lipschitz topology. 
Then, obviously, the sequence $(X, p, \rho_k)$
converges in the pointwise convergence topology
(as defined in 2.3) if and only if 
$\rho_k\in\text{Hom}(\pi,\text{Isom}(X))$
converges in the algebraic topology.

\proclaim{Lemma~2.5} 
Let $\rho_k:\pi\to\text{Isom}(X_k)$
be a sequence of isometric actions
of a discrete group $\pi$ on complete Riemannian
$n$-manifolds $X_k$ such that 
$\rho_k(\pi)$ acts freely.
If the sequence
$(X_k, p_k, \rho_k(\pi))$ converges in the 
equivariant pointed Lipschitz topology to $(X, \Gamma, p)$
and $(X_k, p_k, \rho_k)$ converges to $(X, p, \rho)$ 
in the pointwise convergence topology, then 

$(1)$ $\Gamma$ acts freely, and

$(2)$ $\rho(\pi)\subset\Gamma$, and

$(3)$ $\ker(\rho)\subset\ker(\rho_k)$, for all large $k$.
\endproclaim

\demo{Proof} 
$(1)$ Assume $\gamma\in\Gamma$ and $\gamma(x)=x$.
Choose $\epsilon\in (0,1/10)$ so that there is an $\epsilon$--approximation
$(f_k, g_k, \phi_k, \tau_k)$ of $(X_k, p_k, \rho_k)$ and $(X, p, \Gamma)$
and $x\in B(p,\epsilon/10)$.
Then $g_k(x)=g_k(\gamma(x))=\tau_k(\gamma)(g_k(x))$.
Since $\rho_k(\pi)$ acts freely, $\tau_k(\gamma)=\text{id}$.
By the same argument $\tau_k(\text{id})=\text{id}$.
Hence $\text{id}=\phi_k(\tau_k(\text{id}))=\phi_k(\text{id})=
\phi_k(\tau_k(\gamma))=\gamma$ as desired.

$(2)$ We need to show that $\rho(\gamma)\in\Gamma$,
for any $\gamma\in \pi$. We can assume $\rho(\gamma)\neq\text{id}$.
Choose $\epsilon\in (0,1/10)$ so that the ball $B(p, 1/11\epsilon)$
contains $\rho(\gamma)(p)$ and consider an $\epsilon$--approximation
$(f_k, g_k, \phi_k, \tau_k)$ of $(X_k, p_k, \rho_k)$ and $(X, p, \Gamma)$.

Then for all large enough $k$,
$\rho_k(\gamma)\in B(p_k, 1/10\epsilon)$.
Look at $\tau_k(\rho_k(\gamma))\in\Gamma(1/9\epsilon)$.
Since the set $\Gamma(1/9\epsilon)$ is finite,
we can pass to subsequence to assume that
$\tau_k(\rho_k(\gamma))$ is equal to $\gamma_\epsilon\in\Gamma(1/9\epsilon)$;
thus $\rho_k(\gamma)=\phi_k(\gamma_\epsilon)$.

Take an arbitrary $x\in B(p,1/9\epsilon)$.
Then $g_k(x)\to x$
and, hence, $\rho_k(\gamma)(g_k(x))$ converges to
$\rho(\gamma)(x)$.
Notice that 
$\rho_k(\gamma)(g_k(x))=\phi_k(\gamma_\epsilon)(g_k(x))
\to \gamma_\epsilon(x)$.
So $\rho(\gamma)(x)=\gamma_\epsilon(x)$ 
for any $x\in B(p,1/9\epsilon)$.  

Thus, for any small enough $\epsilon$, we
have found $\gamma_\epsilon\in\Gamma$ that is equal to
$\rho(\gamma)$ on the ball $B(p,1/9\epsilon)$.
Since $\Gamma$ acts freely, 
$\gamma_\epsilon=\gamma_{\epsilon^\prime}$
for all $\epsilon^\prime\le\epsilon$,
that is the element $\gamma_\epsilon\in\Gamma$ is independent
of $\epsilon$. Thus $\rho(\gamma)=\gamma_\epsilon$
everywhere and hence $\rho(\gamma)\in\Gamma$.

$(3)$ Assume $\rho(\gamma)=\text{id}$.
Fix any $\epsilon\in (0,1/10)$. Take $x\in B(p,\epsilon/10)$ and
consider an $\epsilon$--approximation
$(f_k, g_k, \phi_k, \tau_k)$ of $(X_k, p_k, \rho_k)$ and $(X, p, \Gamma)$.
We have $g_k(x)\to x$ and $\rho_k(\gamma)(g_k(x))\to\rho(\gamma)(x)=x$.
Note that $d(x,\phi_k(\rho_k(\gamma))(x))$
is equal to
$$d(f_k(g_k(x)), f_k(\rho_k(\gamma)(g_k(x))))<
e^\epsilon d_k(g_k(x),\rho_k(\gamma)(g_k(x)))@>>k\to\infty >0.$$
Therefore, for all large $k$,
$\phi_k(\rho_k(\gamma))=\text{id}$, because
$\Gamma$ is a discrete subgroup that acts freely. 
Hence $\rho_k(\gamma)=\tau_k(\phi_k(\rho_k(\gamma)))=\tau_k(\text{id})=\text{id}$
as claimed.
\qed\enddemo

\proclaim{Lemma~2.6} 
Let $\rho_k:\pi\to\text{Isom}(X_k)$
be a sequence of isometric actions
of a discrete group $\pi$ on complete Riemannian
$n$-manifolds $X_k$ such that 
$\rho_k(\pi)$ acts freely.
Suppose that the sequence
$(X_k, p_k, \rho_k(\pi))$ converges in the 
equivariant pointed Lipschitz topology to $(X, \Gamma, p)$
and $(X_k, p_k, \rho_k)$ converges to $(X, p, \rho)$ 
in the pointwise convergence topology. 

Then, 
for any $\epsilon>0$ and for any finite subset $S\subset\pi$,
there is $k(\epsilon, S)$ with the property that
for each $k>k(\epsilon, S)$ there exists an
$\epsilon$--Lipschitz approximation 
$(f_k, g_k,\phi_k,\tau_k)$
between $(X_k, p_k, \rho_k(\pi))$
and $(X, p,\Gamma)$ such that 
$\phi_k(\rho_k(\gamma))=\rho(\gamma)$ and 
$\rho_k(\gamma)=\tau_k(\rho(\gamma))$
for every $\gamma\in S$.
\endproclaim

\demo{Proof} Pick $\delta<\epsilon$ so large that
the set $\rho(S)(p)$ is contained in the ball $B(1/10\delta ,p)$.
Choose a $\delta$-Lipschitz approximation
$(f_k, g_k, \phi_k, \tau_k)$ between $(X_k, p_k, \rho_k(\pi))$
and $(X, p,\Gamma)$.
 
Note that for every $\gamma\in S$, the
sequence $\rho_k(\gamma)(g_k(p))$  
converges to $\rho(\gamma)(p)$.
Hence the sequence $f_k(\rho_k(\gamma)(g_k(p)))=
\phi_k(\rho_k(\gamma))(f_k(g_k(p)))=\phi_k(\rho_k(\gamma))(p)$
converges to $\rho(\gamma)(p)$ in $X$.

Since $\Gamma$ is discrete and acts freely,
$\phi_k(\rho_k(\gamma))=\rho(\gamma)$ for large enough $k$.
Moreover, this is true for any $\gamma\in S$,
because $S$ is finite. Also
$\rho_k(\gamma)=(\tau_k\circ\phi_k)(\rho_k(\gamma))=
\tau_k(\rho(\gamma))$.
Thus, this $\delta$-Lipschitz approximation 
has all desired properties.
\qed
\enddemo

\proclaim{Proposition~2.7}
Let $X_k$ be a sequence of Hadamard manifolds
with sectional curvatures in $[a, b]$ for $a\le b<0$
and let $\pi$ be a finitely generated group that 
is not virtually nilpotent.
Let $\rho_k: \pi\to\text{Isom}(X_k)$ be  
an arbitrary sequence of free and isometric actions 
such that $(X_k,p_k,\rho_k)$ converges 
in the pointwise convergence topology.

Then the sequence $(X_k, p_k, \rho_k)$ is precompact
in the equivariant pointed Lipschitz topology.
\endproclaim

\demo{Proof} 
Let $(X, p, \rho)$ be the limit of $(X_k,p_k,\rho_k)$
in the pointwise convergence topology.
Choose $r$ so that the open ball $B(p,r)\subset X$
contains $\{\rho(\gamma_1)(p),\dots\rho(\gamma_m)(p)\}$ where
$\{\gamma_1,\dots\gamma_m\}$ generate $\pi$.
Passing to subsequence, we assume that $B(p_k,r)$
contains $\{\rho_k(\gamma_1)(p),\dots\rho_k(\gamma_m)(p)\}$.

Show that, for every $k$, 
there exists $q_k\in B(p_k,r)$ such that
for any $\gamma\in\pi\setminus\{\text{id}\}$,
we have $\rho_k(\gamma)(q_k)\notin B(q_k,\mu_n/2)$ where
$\mu_n$ is the Margulis constant.
Suppose not. 
Then for some $k$, the whole ball $B(p_k,r)$
projects into
the thin part $\{\text{InjRad}<\mu_n/2\}$ 
under the projection $\pi_k:X_k\to X_k/\rho_k(\pi)$.
Thus the ball $B(p_k, r)$ lies in a connected component $W$
of the $\pi_k$--preimage of the thin part of $X_k/\rho_k(\pi)$. 
According to [2, p111] the stabilizer of
$W$ in $\rho_k(\pi)$ is virtually nilpotent and,
moreover, the stabilizer contains every element
$\gamma\in\rho_k(\pi)$ with $\gamma(W)\cap W\neq\emptyset$.
Therefore, the whole group $\rho_k(\pi)$ stabilizes $W$.
Hence $\rho_k(\pi)$ must be virtually nilpotent.
As $\rho_k$ is injective, $\pi$ is virtually nilpotent.
A contradiction. 
 
Thus, $(X_k, q_k, \rho_k(\pi))$ is 
Lipschitz precompact [20] and, hence passing to subsequence,
one can assume that $(X_k, q_k,\rho_k(\pi))$
converges to some $(X, q, \Gamma)$. 

It is a general fact that follows easily from definitions 
that whenever $(X_k, q_k, \Gamma_k)$ converges to $(X, q, \Gamma)$
in the equivariant pointed Lipschitz topology
and a sequence of points $p_k\in X_k$ converges
to $p\in X$, then
$(X_k, p_k, \Gamma_k)$ converges to $(X, p, \Gamma)$
in the equivariant pointed Lipschitz topology.
\qed\enddemo  

\proclaim{Corollary~2.8}
Let $X_k$ be a sequence of Hadamard manifolds
with sectional curvatures in $[a, b]$ for $a\le b<0$
and let $\pi$ be a finitely generated group that 
is not virtually nilpotent.
Let $\rho_k: \pi\to\text{Isom}(X_k)$ be  
an arbitrary sequence of free and isometric 
such that $(X_k,p_k,\rho_k)$ converges to $(X, p,\rho)$
in the pointwise convergence topology.

Then $\rho$ is a free action, in particular
$\rho$ is injective.
\endproclaim
\demo{Proof}
Pass to a subsequence so that $(X_k, p_k, \rho_k)$
converges to $(X, p, \Gamma)$
in the equivariant pointed Lipschitz topology.
By $2.5(3)$, $\rho$ is injective. 
Furthermore, $\rho(\pi)$ acts freely because it is
a subgroup of $\Gamma$.\qed
\enddemo

\proclaim{Corollary~2.9}
Let $X$ be a complete Riemannian
manifold with sectional
curvatures in $[a, b]$ where $a\le b<0$.
Assume $\pi$ is a torsion free, finitely generated group 
that is not virtually nilpotent.
Consider a subset of faithful discrete 
representations
$S\subset\text{Hom}(\pi,\text{Isom}(X))$
that is precompact in the pointwise convergence
topology. Then, 

$(1)$ for any $p\in X$,
the set $\{(X, p, \rho(\pi)):\rho\in S\}$ 
is precompact in the equivariant pointed Lipschitz topology.

$(2)$ the closure of $S$ in $\text{Hom}(\pi,\text{Isom}(X))$
consists of faithful discrete representations.
\endproclaim
\demo{Proof} Proposition~$2.7$ implies $(1)$ and
Corollary~$2.8$ implies $(2)$. 
\qed\enddemo

\proclaim{Proposition~2.10}
Assume that $\pi$ is a finitely presented discrete group,
that is not virtually nilpotent and does not
have a nontrivial decomposition into an amalgamated product
or an HNN-extension over a virtually nilpotent group.

Let $\rho_k: \pi\to\text{Isom}(X_k)$ be  
an arbitrary sequence of free and isometric actions of $\pi$
on Hadamard $n$-manifolds $X_k$.
Assume that the sectional curvatures of $X_k$
lie in $[a, b]$ for $a\le b<0$.

Then, for some $p_k\in X_k$,the sequence 
$(X_k, p_k, \rho_k)$ is precompact
in the pointwise convergence topology.
\endproclaim

\demo{Proof} Let $S\subset\pi$ be a finite
subset that generates $\pi$ and contains $\{\text{id}\}$.
For $x\in X_k$, we denote $D_k(x)$ the diameter
of the set  $\rho_k(S)(x)$.
Set $D_k=\inf_{x\in X_k}D_k(x)$.

Suppose $D_k$ is unbounded.
Then it follows from a work of Bestvina and Paulin 
[8],[32],[33] (cf.~[27]),
that there exists an action of $\pi$ on a real tree
with no proper invariant subtree 
and virtually nilpotent arc stabilizers.
For completeness we briefly review this construction.
The rescaled pointed Hadamard manifold $\frac{1}{D_k}\!\cdot\! X_k$
has sectional curvature $\le b\cdot D_k\to -\infty$ as $k\to\infty$.
Find $p_k\in X_k$ such that $D_k(p_k)\le D_k+1/k$. 
Consider the sequence of triples 
$(\frac{1}{D_k}\!\cdot\! X_k, p_k, \rho_k)$.
Repeating an argument of Paulin [33,\S 4],
we can pass to subsequence that converges 
to a triple $(X_\infty, p_\infty, \rho_\infty)$.
(For the definition of the convergence see [32],[33].
Paulin calls it ``convergence in the Gromov topology''.) 

The limit space $X_\infty$ is a length space 
of curvature $-\infty$, that is a real tree.
Because of the way we rescaled, 
the limit space has a natural isometric 
action $\rho_\infty$ of $\pi$
with no global fixed point  [32],[33].
Then it is a standard fact that there
exists a unique $\pi$--invariant subtree $T$
of $X_\infty$ that has no proper $\pi$--invariant subtree. 
In fact $T$ is the union of all the axes of all hyperbolic
elements in $\pi$.
Since the sectional curvatures are uniformly bounded
away from zero and $-\infty$, the Margulis lemma implies that
the stabilizer of any 
non-degenerate segment is virtually nilpotent (cf.~[32]). 

Note that any increasing sequence of virtually
nilpotent subgroups of $\pi$ is stationary.
Indeed, since a virtually nilpotent group is amenable, 
the union $U$ of an increasing sequence 
$U_1\subset U_2\subset U_3\subset\dots$
of virtually nilpotent subgroups is also an amenable group.
If the fundamental group of a complete
manifold of pinched negative curvature is amenable,
it must be finitely generated [13], [10].
In particular, $U$ is finitely generated,
hence $U_n=U$ for some $n$. 
Thus, the $\pi$--action on the tree $T$ is stable
[9, Proposition 3.2(2)].

We summarize that the $\pi$--action
on $T$ is stable, has virtually nilpotent arc stabilizers
and no proper $\pi$--invariant subtree.
Therefore, the Rips machine [9, Theorem 9.5] 
produces a splitting of
$\pi$ over a virtually solvable group.
Any amenable subgroup of $\pi$ is virtually nilpotent [13], [10],
hence $\pi$ splits over a virtually nilpotent group.
This is a contradiction with the assumption that
$D_k$ is unbounded.

Thus the sequence $D_k(p_k)$ is bounded, therefore
as we observed in $2.3$, the sequence $(X_k, p_k,\rho_k)$ 
is precompact in the pointwise convergence topology.
\qed\enddemo

\head 
\S 3.~Vector bundles and $\widetilde{KO}$-theory.
\endhead

This section provides some $\widetilde{KO}$-theoretic background. 
All of the facts below are well-known to experts; 
however it is usually not easy to locate a reference.

In this paper we deal with 
rank $n$ real vector bundles over open smooth $n$--manifolds.
It is well-known that any open smooth manifold $N$ is homotopy equivalent to
CW--complex of dimension $<\dim(N)$, hence
two rank $n$ vector bundles over $N$ are isomorphic
iff they are stably isomorphic [26, 8.1.5].
The set of stable isomorphism classes of real vector bundles
over a space $X$ forms an abelian group $\widetilde{KO}(X)$.
The correspondence $X\to\widetilde{KO}(X)$ is clearly
a functor. This functor is isomorphic to the functor $[-,BO]$
on the category of connected CW-complexes of uniformly bounded 
dimension [26, 8.4.2]. 

We now recall the definition of the $\widetilde{KO}^*$--theory. 
By the Bott periodicity, 
$\Omega^8(BO\times\Bbb Z)$ is weakly homotopy
equivalent to $BO\times\Bbb Z$ [37, 11.60].                                                                                                                                                               
We can use the spectrum $\Omega^n(BO\times\Bbb Z)$
to define a generalized (reduced) cohomology theory 
$\widetilde{KO}^*$ on the category of all pointed
connected finite--dimensional CW-complexes and point-preserving maps.
In other words,
we set $\widetilde{KO}^n(X,x)=[X,x;\Omega^n(BO\times\Bbb Z),\ast]$ 
(cf.~[37, 8.42]).

Let $F$ be a forgetful map from the category of pointed
connected CW--complexes to the category of connected CW--complexes.
We now show that the functors $\widetilde{KO}^0(-)$ 
and $\widetilde{KO}(F(-))$ are isomorphic on 
the category of pointed connected CW-complex\-es of uniformly bounded 
dimensions.
Indeed, since we deal with connected CW-complexes, the functor
$\widetilde{KO}^0(-)=[-;BO\times\Bbb Z),\ast]$ is isomorphic
to the functor $[-;BO,\ast ]$.
As we mentioned above $\widetilde{KO}(-)\cong [-,BO]$.
Therefore, it remains to show that the transformation
of functors
$[-;BO,\ast ]\to [F(-);BO]$ is an isomorphism.
In other words, it suffices to check that
the map $[X,x;BO,\ast ]\to [F(X,x);BO]=[X,BO]$  
of based homotopy classes into free homotopy
classes is bijective. 
Indeed, any map $X\to BO$ is homotopic to
a basepoint preserving map [36, 7.3.Lemma~2]
which is unique up to based homotopy 
because $BO$ is a path-connected H-space. 
[36, 7.3.Theorem~5]. 
Thus, $\widetilde{KO}^0(-)\cong\widetilde{KO}(F(-))$.
Most of the time we suppress the base points
and treat the functors $\widetilde{KO}^0(-)$ and
$\widetilde{KO}(-)$ as isomorphic.

Rationally the group
$\widetilde{KO}^i(X)$ can be easily computed
in terms of cohomology of $X$ as explained 
in the lemma below.
Recall that the Bott periodicity implies that,
after tensoring with rationals,
the $\widetilde{KO}$--groups
of the $0$--sphere are given by
$\widetilde{KO}^i(S^0)\otimes\Bbb Q\cong\Bbb Q$
if $i\equiv 0\pmod{4}$ 
and $\widetilde{KO}^i(S^0)\otimes\Bbb Q=0$ otherwise [26, 15.12.3]. 

\proclaim{Lemma~3.1} 
On the category of pointed connected CW-complexes of 
uniformly bounded dimension 
the contravariant functors $\widetilde{KO}^i(-)\otimes\Bbb Q$ and 
$\oplus_{n>0}H^{n}(-)\otimes
\widetilde{KO}^{i-n}(S^0)\otimes\Bbb Q$ are isomorphic.
\endproclaim

\demo{Proof} 
Consider the generalized cohomology
theory $\widetilde{H}^*(-)\otimes \widetilde{KO}^*(S^0)\otimes\Bbb Q$
where $\widetilde{H}^*(-)$ is the ordinary reduced cohomology.
Its coefficients are isomorphic to the coefficients of the theory
$\widetilde{KO}^*(-)\otimes\Bbb Q$. 
Since the coefficients are rational vector spaces,
[24, 3.22(ii)] implies that the theories are naturally
equivalent on the category of pointed finite-dimensional
CW--complexes. 

For any connected finite-dimensional complex $X$,
$\widetilde{H}^n(X)=0$, for all $n\le 0$.
Clearly for $n>0$, the functors
$\widetilde{H}^n(-)$ and $H^n(-)$
are isomorphic. 
Hence we have an isomorphism of functors
$$\widetilde{KO}^i(-)\otimes\Bbb Q\cong
\oplus_{n\ge 0}\widetilde{H}^{n}(-)\otimes
\widetilde{KO}^{i-n}(S^0)\otimes\Bbb Q
\cong\oplus_{n>0}H^n(-)\otimes\widetilde{KO}^{i-n}(S^0)\otimes\Bbb Q.\qed$$
\enddemo

\subheading{Remark~3.2} In particular,
we have isomorphisms of functors
$$\widetilde{KO}^0(-)\otimes\Bbb Q\cong
\oplus_{n>0}H^{4n}(-,\Bbb Q)\ \ \text{and}\ \ 
\widetilde{KO}^{-1}(-)\otimes\Bbb Q\cong
\oplus_{n>0}H^{4n-1}(-,\Bbb Q).$$

\proclaim{Corollary~3.3}
Let $N$ be an open smooth manifold
such that the torsion subgroup of
$\widetilde{KO}(N)$ is finite
and let $n\ge\dim (N)$.

Then the set of isomorphism classes of
rank $n$ vector bundles over $N$ is infinite 
if and only if $H^{4k}(N,\Bbb Q)\neq 0$ for some $k>0$.
\endproclaim

\demo{Proof} Any open manifold $N$ is homotopy equivalent to
a CW-complex of dimension $<\dim(N)$,
therefore, two rank $n$ vector bundles over $N$ are
isomorphic iff they are stably isomorphic. The result now
follows from the above remark.\qed
\enddemo

\proclaim{Lemma~3.4}
Let $Y$ be a connected finite-dimensional
CW-complex such that $\pi_1(Y)$ is finitely generated and
$\dim_{\Bbb Q}H^{4k}(Y,\Bbb Q)<\infty$ for all $k>0$.
Then there exists a finite connected subcomplex $X\subset Y$
such that the inclusion $i:X\to Y$ induces
a surjection on the fundamental groups
and an injection 
$\widetilde{KO}(Y)\otimes\Bbb Q\to\widetilde{KO}(X)\otimes\Bbb Q$
\endproclaim

\demo{Proof} 
It is a standard fact [19, p117] that for any homology
class $\alpha\in H_i(Y,\Bbb Q)$ there exists a finite
(not necessarily connected)
CW-complex $X_\alpha$ and a continuous map 
$f_\alpha:X_\alpha\to Y$ such that 
$\alpha\in f_{\alpha*}(H_i(X_\alpha, \Bbb Q))$.
For every $k>0$, choose a finite basis
in the $\Bbb Q$--vector space $H_{4k}(Y,\Bbb Q)$
and construct such a finite CW--complex $X_\alpha$
for every basis element $\alpha$.

Let $X_0$ be 
the the disjoint union of all these complexes over all $k$
and all the basis elements.
Since $Y$ is finite-dimensional with finite $4k$-th Betti numbers,
the CW--complex $X_0$ is finite.
Note that $X_0$ comes with a 
continuous map $f:X_0\to Y$ that induces a surjection
on rational $4k$-th homology, and therefore,
an injection on rational $4k$-cohomology for all $k>0$.

Finally, set $X$ to be an arbitrary
connected finite subcomplex of $Y$ such that
$f(X_0)\subset X$.
Using that $\pi_1(Y)$ is finitely generated,
we add to $X$ finitely many $1$-cells of $Y$ to make sure
the inclusion $i:X\hookrightarrow Y$ 
induces a surjection of fundamental groups.
Moreover, by construction
$i$ induces an injection on rational $4k$-cohomology for all $k>0$.
By the remark $3.2$,  
$\widetilde{KO}(Y)\otimes\Bbb Q\to\widetilde{KO}(X)\otimes\Bbb Q$
is injective.
\qed\enddemo

\head
\S 4.~An invariant of actions.
\endhead

Let $K$ be a finite-dimensional connected CW-complex with a 
reference point $q$. 
Denote $\tilde{K}\to K$ the universal cover of $K$;
choose $\tilde{q}\in \tilde{K}$ that is mapped to $q$
by the covering projection.
Consider a pointed contractible manifold $X$
and an arbitrary action 
$\rho:\pi_1(K,q)\to\text{Diffeo}(X)$.

Since $X$ is contractible, the $X$-bundle 
$\tilde{K}\times_{\rho} X$ over $K$ has a section.
Any two sections are homotopic through sections. 

We now define a certain invariant of actions of 
$\pi_1(K,q)$ on $X$. Given such an action $\rho$,
consider the vertical vector bundle 
$\tilde{K}\times_\rho TX$ over $\tilde{K}\times_\rho X$
and set $\tau(\rho)$ to be the pullback of the vertical
bundle via an arbitrary section $s:K\to \tilde{K}\times_\rho X$
Clearly, two actions that are conjugate 
in $\text{Diffeo}(X)$ have same invariants. 

Any section can be lifted to
a $\rho$-equivariant continuous map 
$\tilde{K}\to \tilde{K}\times X\to X$.
Any two $\rho$-equivariant continuous maps 
$\tilde{g}, \tilde{f}:\tilde{K}\to X$, 
are $\rho$-equivariantly homotopic.
Indeed, $\tilde{f}$ and $\tilde{g}$ descend to sections
$K\to \tilde{K}\times_\rho X$ that must be homotopic.
This homotopy lifts to a $\rho$--equivariant homotopy
of $\tilde{f}$ and $\tilde{g}$. 
 
Assume now that $\rho(\pi_1(K,q))$ acts freely
and properly discontinuously on $X$,
so the map $\pi :X\to X/\rho(\pi_1(K))$ is a covering.
Then the map $\tilde{f}$ descends to a continuous map 
$f:K\to X/\rho(\pi_1(K))$.
Thus, to any action $\rho:\pi_1(K,q)\to\text{Homeo}(X)$
such that $\rho(\pi_1(K,q))$ acts freely
and properly discontinuously on 
a contractible manifold $X$, we associate a 
(unique up to homotopy) continuous map $f$.
Observe that if $\rho$ is injective,
then $f$ induces an isomorphism of fundamental groups.

We now observe that, 
in case $\rho(\pi_1(K,q))$ acts freely
and properly discontinuously on $X$, 
the invariant $\tau(\rho)$ is equal to $f^*TX/\rho(\pi_1(K))$,
the pullback of the tangent bundle to $X/\rho(\pi_1(K))$
via $f$. 

If we consider only orientation-preserving actions
the vector bundle $\tau(\rho)$ comes with a natural 
orientation. 

\head
\S 5.~Main results
\endhead

Throughout this section $K$ is a finite-dimensional, 
connected CW-complex with a reference point $q$.
Let $\tilde{K}$ be the universal cover of $K$,
$\tilde{q}\in\tilde{K}$ be a preimage of $q\in K$. 
Using the point $\tilde{q}$ we identify $\pi_1(K,q)$ with
the group of automorphisms of the covering $\tilde{K}\to K$.
Recall that $\tau$ is an invariant of actions
defined in $\S 4$.

\proclaim{Theorem~5.1} Let the complex $K$ is finite.
Let $\rho_k:\pi=\pi_1(K, q)\to\text{Isom}(X_k)$
be a sequence of isometric actions
of $\pi_1(K)$ on Hadamard $n$-manifolds $X_k$ such that
$\rho_k(\pi)$ is a discrete subgroup of $\text{Isom}(X_k)$
that acts freely. 

Suppose that, for some $p_k\in X_k$,
$(X_k, p_k, \rho_k(\pi))$ converges in the
equivariant pointed
Lipschitz topology to $(X, p, \Gamma)$
and $(X_k, p_k, \rho_k)$ converges to $(X, p, \rho)$
in the pointwise convergence topology. 

Then $\tau(\rho_k)=\tau(\rho)$ for all large $k$.
\endproclaim 

\demo{Proof} 
Let $\tilde{q}\in F\subset\tilde{K}$ be a finite subcomplex
that projects onto $K$.
Clearly, the finite set 
$S=\{\gamma\in\pi_1(K,q): \gamma(F)\cap F\neq\emptyset\}$
generates $\pi_1(K,q)$.

Note that $X$ is contractible, indeed
any spheroid in $X$ lies in the diffeomorphic
image of a metric ball in $X_j$.
Any metric ball in a Hadamard manifold is contractible.
Thus $\pi_*(X)=1$. 

Using $\S 4$, we find a $\rho$-equivariant continuous map
$\tilde{h}: \tilde{K}\to X$.
Recall that by definition $\tau(\rho)$ is the pullback via 
$h:K\to X/\rho(\pi_1(K,q)$
of the tangent bundle to $X/\rho(\pi_1(K,q)$. 
Let $\bar{h}:K\to X/\Gamma$ be the composition of
$h$ and the covering $X/\rho(\pi_1(K,q)\to X/\Gamma$.
Therefore, $\tau(\rho)$ is also the pullback
via $\bar{h}$ of the tangent bundle to $X/\Gamma$.  

Choose $\epsilon >0$ so small that $\tilde{h}(F)$ lies in the open
ball $B(p,1/10\epsilon)\subset X$.
For large $k$, we find an $\epsilon$-Lipschitz
approximation $(\tilde{f}_k, \tilde{g}_k, \phi_k, \tau_k)$ 
between $(X_k,p_k,\rho_k(\pi))$ and
$(X,p,\Gamma)$.  By lemma~$2.6$ we can assume that
$\tau_k(\rho(\gamma ))=\rho_k(\gamma)$ for all $\gamma\in S$.
Hence, the map
$\tilde{h}_k=\tilde{g}_k\circ\tilde{h}:F\to X_k$
is $\rho_k$-equivariant.
Extend it by equivariance to a $\rho_k$-equivariant
map $\tilde{h}_k:\tilde{K}\to X_k$.
Passing to quotients we get a map
$h_k:K\to X_k/\rho_k(\pi_1(K))$
such that $\tau(\rho_k)$ is the pullback via $h_k$
of the tangent bundle to $X_k/\rho_k(\pi_1(K))$.

By construction, $h_k=g_k\circ \bar{h}$
where $g_k$ is the drop of $\tilde{g}_k$.
Being a diffeomorphism $g_k$ preserves tangent bundles.
Therefore, the pullback via $h_k$
of the tangent bundle to $X_k/\rho_k(\pi_1(K))$
is equal to the pullback via $\bar{h}$
of the tangent bundle to $X/\rho(\pi_1(K,k)$.
In other words $\tau(\rho_k)=\tau(\rho)$.
The proof is complete.\qed
\enddemo

\subheading{Remark~5.2} In the above theorem it is possible to
keep track of orientations provided 
all the actions $\rho_k$
on $X_k$ preserve orientations  
(it makes sense because
being a contractible manifold $X_k$ is orientable).
Indeed, fix an orientation on $X$ (which is also
contractible) and choose orientations of $X_k$
so that diffeomorphisms $g_k$ preserve orientations.
Then, obviously, the theorem~$5.1$ is still true
where the vector bundle isomorphism $\tau(\rho_k)=\tau(\rho)$
preserves orientations.

\proclaim{Theorem~5.3} 
Assume that 
$\dim_{\Bbb Q}\widetilde{KO}(K)\otimes\Bbb Q<\infty$.
Let $\pi_1(K)$ be a finitely generated group and
$\rho_k:\pi=\pi_1(K)\to\text{Isom}(X_k)$
be a sequence of isometric actions
of $\pi$ on Hadamard $n$-manifolds $X_k$ such that
$\rho_k(\pi)$ is a discrete subgroup of $\text{Isom}(X_k)$
that acts freely. 

Suppose that, for some $p_k\in X_k$,
$(X_k, p_k, \rho_k(\pi))$ converges in the
equivariant pointed
Lipschitz topology to $(X, p, \Gamma)$
and $(X_k, p_k, \rho_k)$ converges to $(X, p, \rho)$
in the pointwise convergence topology. 

Then, for all large $k$, the images of $\tau(\rho_k)$
and $\tau(\rho)$ in $\widetilde{KO}(K)\otimes\Bbb Q$
are equal.
\endproclaim
\demo{Proof} 
By $3.4$ there exists a finite connected subcomplex
$X\subset K$ such that the inclusion $i:X\to K$
induces a $\pi_1$--epimorphism
and a $\widetilde{KO}$-monomorphism.

The sequence of isometric actions 
$\rho_k\circ i_*$ of $\pi_1(X,x)$ on $X_k$ satisfies
the assumptions of the theorem~$5.1$, therefore,
$\tau(\rho_k\circ i_*)=\tau(\rho\circ i_*)$ for all large $k$.
In other words, the pullback bundles $i^*\tau(\rho_k)$
are isomorphic to the vector bundle $i^*\tau(\rho)$.
Since $i$ induces a monomorphism of rational $\widetilde{KO}$--groups,
the images of $\tau(\rho_k)$
and $\tau(\rho)$ in $\widetilde{KO}(K)\otimes\Bbb Q$
are equal for all large $k$.
\qed\enddemo

\proclaim{Corollary~5.4}
Assume that the group
$\widetilde{KO}(K)$ is finitely generated.
Let $\pi_1(K)$ be a finitely generated group
and let $\rho_k:\pi=\pi_1(K)\to\text{Isom}(X_k)$
be a sequence of isometric actions
on Hadamard $n$-manifolds $X_k$ such that 
$\rho_k(\pi)$ is a discrete subgroup of $\text{Isom}(X_k)$
that acts freely. 
Suppose that, for some $p_k\in X_k$,
$(X_k, p_k, \rho_k)$ is precompact in both
pointwise convergence topology and equivariant pointed
Lipschitz topology. 

Then the set $\{\tau(\rho_k)\}$ falls into finitely many
stable isomorphism isomorphism classes.
\endproclaim 

\demo{Proof} 
Argue by contradiction. 
Pass to subsequence to assume that
all $\tau(\rho_k)$ are different and that
$(X_k, p_k, \rho_k(\pi_1(K,q)))$ converges to $(X, p, \Gamma)$
in the equivariant pointed Lipschitz topology and
$(X_k, p_k, \rho_k)$
converges to $(X, p, \rho)$ in the pointwise convergence 
topology.

Note that if the group
$\widetilde{KO}(K)$ is finitely generated, then
$\dim_{\Bbb Q}\widetilde{KO}(K)\otimes\Bbb Q)<\infty$.
Then by the theorem~$5.3$ the images of $\tau(\rho_k)$ in 
$\widetilde{KO}(K)\otimes\Bbb Q$
are all equal for large $k$. 
In other words for large $k$, the images of $\tau(\rho_k)$ in 
$\widetilde{KO}(K)$ are all equal modulo torsion.
Since the abelian group, 
$\widetilde{KO}(K)$ is finitely generated,
its torsion subgroup is finite. 
Therefore, passing to subsequence we can assume
that the images $\tau(\rho_k)$ in 
$\widetilde{KO}(K)$ are the same, 
a contradiction.\qed
\enddemo

\proclaim{Corollary~5.5}
Assume that the complex $K$ is aspherical.
Suppose that the groups $\widetilde{KO}(K)$ and $\pi_1(K)$ are 
finitely generated.
Let $\rho_k:\pi=\pi_1(K, q)\to\text{Isom}(X_k)$
be a sequence of free, isometric actions
of $\pi_1(K)$ on Hadamard $n$-manifolds $X_k$. 
Suppose that, for some $p_k\in X_k$,
$(X_k, p_k, \rho_k(\pi))$ is precompact in the
equivariant pointed Lipschitz topology and 
$(X_k, p_k, \rho_k)$  is precompact in 
the pointwise convergence topology. 

Then the set of manifolds $\{X_k/\rho_k(\pi)\}$
falls into finitely many tangential equivalence types.
\endproclaim

\demo{Proof}
By the corollary~$5.4$ the set 
$\{\tau(\rho_k)\}$ falls into finitely many
stable isomorphism isomorphism classes.
Let $\tilde{h_k}:\tilde{K}\to X_k$ be $\rho_k$--equivariant
maps constructed in $\S 3$. 
Passing to quotients we get homotopy equivalences
$h_k:K\to X_k/\rho_k(\pi)=N_k$
with the property that $h_k^*TN_k\cong\tau(\rho_k)$.
Therefore, $X_k/\rho_k$ and $X_j/\rho_j$ are tangential homotopy 
equivalent iff $\tau(\rho_k)\cong\tau(\rho_j)$.
\qed\enddemo

\proclaim{Corollary~5.6}
Suppose that the complex $K$ is aspherical and
the group $\widetilde{KO}(K)$ is finitely generated.
Let $\pi_1(K)$ be a finitely presented group 
that is not virtually nilpotent.
Assume that $\pi_1(K)$ does not split 
as a nontrivial amalgamated product or
an HNN-extension 
over a virtually nilpotent group.

Then, for any $a\le b<0$ and an integer $n\ge 2$,
the class $\Cal M_{a, b, \pi_1(K), n}$ breaks into finitely
many tangential homotopy types.
\endproclaim

\demo{Proof} Apply the corollary~$5.5$ and the proposition $2.7$, $2.10$.\qed
\enddemo

\proclaim{Corollary~5.7} Let $\pi$ be the fundamental group
of a finite aspherical CW-complex.
Suppose that $\pi$ is not virtually nilpotent and that $\pi$
does not split as a nontrivial amalgamated product or
an HNN-extension 
over a virtually nilpotent group.

Then, for any $a\le b<0$ and an integer $n\ge 2$,
the class $\Cal M_{a, b, \pi_1(K), n}$ breaks into finitely
many tangential homotopy types.
\endproclaim

\demo{Proof} Let $\pi$ be the
the fundamental group of a finite complex $K$.
Then $\pi$ is finitely presented and the group $\widetilde{KO}(K)$
is finitely generated [24, p52]. 
Hence the corollary~$5.6$ applies.
\qed\enddemo

\head
\S 6.~Graphs of groups and Accessibility.
\endhead

In this section we explain how to generalize
Theorem~$1.1$ to certain amalgamated products and
HNN--extensions, or more generally,
to some graphs of groups.
Then we deduce the corollary~$1.2$.

Recall that a graph of groups is a graph
whose vertices and edges are labeled
with {\it vertex groups} $\pi_v$ and 
{\it edge groups} $\pi_e$ and such that
every pair $(v,e)$ where 
the edge $e$ is incident to the vertex $v$
is labeled with a group monomorphism
$\pi_e\to\pi_v$.
We only consider finite connected graphs of groups.
To each graph of groups one can associate its 
fundamental group which is a result of repeated
amalgamated products and HNN--extensions of vertex groups
over the edge groups (see [3] for more details). 

\proclaim{Proposition~6.1}
Let $\pi$ be the fundamental group of a finite graph of 
groups with the vertex 
groups $\pi_v$.
Assume that the homomorphism
$$\widetilde{KO}(K(\pi,1))\to
\oplus_v\widetilde{KO}(K(\pi_v,1))$$
induced by the inclusions $\pi_v\to\pi$ has finite kernel.

Fix $n$ and $a\le b<0$ and suppose that, 
for each vertex group $\pi_v$,
the class $\Cal M_{a, b,\pi_v, n}$ breaks into finitely
many tangential homotopy types. 
Then so does the class $\Cal M_{a, b, \pi, n}$. 
\endproclaim

\demo{Proof}
Assume, by contradiction that there exist a sequence 
$N_k\in\Cal M_{a, b,\pi, n}$ of manifolds that are not
pairwise tangentially homotopy equivalent. 
It defines an infinite sequence of distinct elements in
$\widetilde{KO}(K(\pi,1))$.
 
The homomorphism 
$\widetilde{KO}(K(\pi,1))\to
\oplus_v\widetilde{KO}(K(\pi_v,1))$
has finite kernel and the
set of vertices is finite, 
therefore, for some vertex $v$ we get an infinite
sequence of elements in 
$\widetilde{KO}(K(\pi_v,1))$.
Each of elements of the infinite
sequence comes from the tangent bundle
to a manifold $\bar{N}_k\in\Cal M_{a, b, \pi_v, n}$,
namely, $\bar{N}_k$ is the covering of $N_k$
induced by the inclusion $\pi_v\to\pi$.
Thus, for this vertex $\Cal M_{a, b, \pi_v, n}$
falls into infinitely many tangential homotopy types,
a contradiction.\qed
\enddemo

\proclaim{Corollary~6.2}
Let $\pi$ be the fundamental group of a finite graph of groups
such that each edge group is a virtually nilpotent
group of cohomological dimension $\le 2$.

Fix $n$ and $a\le b<0$ and suppose that, 
for each vertex group $\pi_v$,
the class $\Cal M_{a, b, \pi_v, n}$ breaks into finitely
many tangential homotopy types. 
Then so does the class $\Cal M_{a, b, \pi, n}$. 
\endproclaim

\demo{Proof}
We can assume that $\Cal M_{a, b, \pi, n}\neq\emptyset$.
Then $\Cal M_{a, b,\pi_e, n}\neq\emptyset$
for every edge group $\pi_e$.
Hence all the edge groups are finitely generated [10].
Being a finitely generated virtually nilpotent group, 
each edge group $\pi_e$ is a fundamental group
of a closed aspherical manifold [17]
which has to be of dimension $\le 2$ 
since $\text{cd}(\pi_e)\le 2$.
 
In particular, the group $\widetilde{KO}^{-1}(K(\pi_e,1))$
is a finitely generated abelian group [24, p52].
According to $3.2$, there is an isomorphism of functors
$\widetilde{KO}^{-1}(-)\otimes\Bbb Q
\cong\oplus_{k>0}H^{4k-1}(-,\Bbb Q)$
on the category of connected CW--complexes of uniformly bounded
dimension. 
Hence, for each edge group $\pi_e$,
the group $\widetilde{KO}^{-1}(K(\pi_e,1))$
is finite.

Following [34] we can assemble the cell complexes
$K(\pi_v,1)$ and $K(\pi_e,1)\times [-1,1]$
into an $K(\pi,1)$ cell complex by using edge--to--vertex
monomorphisms.
Consider the $\widetilde{KO}^*$--Mayer-Vietoris
sequence applied to the $K(\pi,1)$ complex (cf.~[12, VII.9]).
The group 
$\oplus_{e}\widetilde{KO}^{-1}(K(\pi_e,1))$
is finite, so by exactness,
$\widetilde{KO}(K(\pi,1))\to
\oplus_v\widetilde{KO}(K(\pi_v,1))$
has finite kernel.
We now apply $6.1$ to complete the proof.
\qed
\enddemo

\subheading{6.3.~An accessibility result} 
Delzant and Potyagailo have recently proved a powerful
accessibility result which we state below for reader's
convenience.
The definition~$6.4$ and the theorem~$6.5$ 
are taken from [16].

\subheading{Definition~6.4}
A class $\Cal E$ of 
subgroups of a group $\pi$
is called {\it elementary} provided the following four
conditions hold.

(i) $\Cal E$ is closed under conjugation in $\pi$;

(ii) any infinite group from $\Cal E$ is contained
in a {\it unique} maximal subgroup from $\Cal E$;
  
(iii) if a group from the class $\Cal E$ acts on a tree,
it fixes a point, an end, or a pair of ends;

(iv) each maximal subgroup from $\Cal E$
is equal to its normalizer in $\pi$.

Note that the condition (iii) holds
if any group in $\Cal E$ is amenable [31]. 

\proclaim{Theorem~6.5\ \rm [16]\it} 
Let $\pi$ be a finitely presented group without $2$-torsion
and let $\Cal E$ be an elementary class of subgroups of $\pi$.
Then there exists an integer $K>0$ and a {\bf finite} sequence 
$\pi_0,\pi_1,\dots ,\pi_m$ of subgroups of $\pi$ such that

$(1)$ $\pi_m=\pi$, and

$(2)$ for each $k$ with $0\le k<K$,
the group $\pi_k$ either belongs to $\Cal E$ or
does not split as a nontrivial
amalgamated product or an HNN-extension over
a group from $\Cal E$, and

$(3)$ for each $k$ with $K\le k\le m$, the group $\pi_k$ 
is the fundamental group of a finite graph of groups
with edge groups from $\Cal E$, vertex groups from
$\{\pi_0,\pi_1,\dots ,\pi_{k-1}\}$,
and proper edge-to-vertex homomorphisms.
\endproclaim

\proclaim{Proposition~6.6} 
If $\Cal M_{a, b, \pi, n}\neq\emptyset$, then 
the class of virtually nilpotent subgroup of $\pi$
is elementary.
\endproclaim
\demo{Proof} The conditions (i)  
is trivially satisfied.
Any virtually nilpotent group is amenable, hence
(iii) holds thanks to [31].

Let $X$ be a Hadamard manifold such that
$X/\pi\in \Cal M_{a, b, \pi, n}$.
Being the fundamental group of an aspherical manifold,
the group $\pi$ is torsion free.
According to [11], a nontrivial torsion-free virtually nilpotent
discrete isometry group $\Gamma$ of $X$ is characterized
by its fixed-point-set at infinity. 
In fact, either the fixed-point-set is a point
or it consists of two points.
Conversely, any discrete torsion-free subgroup of 
$\text{Isom}(X)$ that fixes a point 
at infinity is virtually nilpotent.
Thus, a virtually nilpotent subgroup $N$ of $\pi$
is maximal if and only if $N$ contains every element of $\pi$
that fixes the fixed-point-set of $N$.
This characterization proves (ii).

Finally, we deduce (iv). Let $g$ be an element 
of the normalizer of a maximal virtually nilpotent
subgroup $N$. Then $g$ preserves the fixed-point-set of $N$
setwise. If the fixed-point-set is a point,
we conclude $g\in N$.
If the fixed-point-set consists of two points $p$ and $q$
and $g\notin N$,
then $g(p)=q$.
Hence $g$ must fix a point on the geodesic
that joins $p$ and $q$, so the element $g$ is elliptic.
This is a contradiction because no discrete torsion-free
group contains an elliptic element.
\qed
\enddemo

\proclaim{Theorem~6.7}
Let $\pi$ be a finitely presented group
such that any 
nilpotent subgroup of $\pi$ has cohomological dimension $\le 2$.
Assume that the group 
$\widetilde{KO}(K(\pi,1))$ is finitely generated.

Then, for any positive integer $n$ and any negative reals
$a\le b$, the class $\Cal M_{a, b, \pi, n}$ 
breaks into finitely
many tangential homotopy types.
\endproclaim

\demo{Proof}
Applying the theorem~$6.5$, 
we get a sequence $\pi_0,\pi_1,\dots ,\pi_m$ 
of subgroups of $\pi$.
In particular, for every $k$ with $0\le k<K$, the group $\pi_k$ 
does not split as a nontrivial
amalgamated product or an HNN-extension over
a virtually nilpotent subgroup of $\pi$.
By the proposition~$6.9$ below, the group
$\pi_k$ is finitely presented and  
$\widetilde{KO}(K(\pi,1))$ is finitely generated.
Therefore, by the theorem~$5.6$,
the class $\Cal M_{a, b, \pi_k, n}$ breaks into finitely
many tangential homotopy types.

Repeatedly applying the corollary~$6.2$, 
we deduce that $\Cal M_{a, b, \pi, n}$ 
breaks into finitely
many tangential homotopy types.
\qed\enddemo

\proclaim{Corollary~6.8}
Let $\pi$ be a word-hyperbolic group.
Then for any $n$ and $a\le b<0$, 
the class $\Cal M_{a, b, \pi, n}$ breaks into finitely
many tangential homotopy types.  
\endproclaim

\demo{Proof}
We can assume that $\Cal M_{a, b, K(\pi, 1), n}\neq\emptyset$,
hence $\pi$ is torsion--free.
Any torsion--free word--hyperbolic group is 
the fundamental group of
a finite CW-complex $K$[15, 5.24], 
in particular $\pi$ is finitely presented and
has finitely generated $\widetilde{KO}(K)$ [24, p52]. 
Any nilpotent subgroup of a torsion free word hyperbolic group
is either trivial or infinite cyclic (see e.g.[4]).
The result now follows from the previous theorem.
\qed\enddemo

\proclaim{Proposition~6.9}
Let $\pi$ be the fundamental group of a finite 
graph of groups with virtually nilpotent edge groups.
Assume $\Cal M_{a, b, \pi, n}\neq\emptyset$. Then

$(1)$ $\pi$ is finitely presented iff
all the vertex groups are finitely presented, and

$(2)$ $\dim_{\Bbb Q}\oplus_{k>0} H^{4k}(\pi, \Bbb Q)<\infty$ iff
$\dim_{\Bbb Q}\oplus_{k>0} H^{4k}(\pi_v, \Bbb Q)<\infty$
for every vertex group $\pi_v$.

$(3)$ the group $\widetilde{KO}(\pi)$ is finitely generated
iff for every vertex group $\pi_v$
the group $\widetilde{KO}(\pi)$ is finitely generated.
\endproclaim
\demo{Proof}
All the edge groups are finitely generated [10].
Being a finitely generated virtually nilpotent group, 
each edge group $\pi_e$ is a fundamental group
of a closed aspherical manifold [17].
In particular, $\widetilde{KO}^{*}(\pi_e)$ as well as
$H^*(\pi_e)$
is a finitely generated abelian group.
Hence, the parts $(2)$ and $(3)$ follow from the
Mayer-Vietoris sequence (cf.~[12, VII.9]).

We now prove $(1)$.
It trivially follows from definitions that if all 
vertex groups are finitely presented,
then so is $\pi$.

Assume now that $\pi$ is finitely presented.
Then every vertex group $\pi_v$ 
finitely generated [3, Lemma~13,p.158].
The fundamental
group $\pi$ of a graph of groups has a  
presentation that can be described as follows.
The set of generators $S$ for for this presentation
$\pi$ is the union of the (finite) sets of generators 
of all the vertex groups; so $S$ is finite.
The set of relations $R$ is a union of the sets of relations 
of all the vertex groups and (finitely many) relations coming
from amalgamations and HNN-extensions over the edge groups.

Since $\pi$ is finitely presented, the group $\pi$
has a presentation $\langle S|R^\prime\rangle$ on the same set of 
generators $S$ where $R^\prime$ is a finite subset of $R$
[3, Theorem~12, p52].
By adding finitely many (redundant) relations, we can assume that
all the relations coming from the edge groups 
still belong to the new set of relations $R^\prime$.
This defines a new graph of groups decomposition of $\pi$
with the same underlying graph, same edge groups and
with finitely presented vertex groups $\pi_v^\prime$.
By construction, the subgroups $\pi_v^\prime$
and $\pi_v$ of $\pi$ have the same set of generators,
hence $\pi_v^\prime=\pi_v$.
Thus, we have found a finite 
presentation for every vertex group.
\qed\enddemo

\head
\S 7.~Finiteness for thickenings
\endhead

Throughout this section $K$ is a finite, 
connected CW-complex. 

By an $n$--{\it thickening} of $K$ we mean 
a compact smooth manifold $L$ of dimension 
$\ge\dim(K)+3$ such that $L$ is simply homotopy 
equivalent to $K$, and 
the inclusion $\partial L\to L$ induces
an isomorphism of fundamental groups (see [38]).
To avoid low-dimensional complications, 
we always assume that $n>4$.

A thickening $L$ of $K$ is always diffeomorphic to the 
regular neighborhood
of a finite simplicial subcomplex $K^\prime\subset L$
with $\dim(K^\prime)\le\dim(K)$ [25, p219][35].
 
If $\pi_1(K)$ is isomorphic
to the fundamental group of a
complete manifold of sectional curvature pinched
between two negative constants,
any homotopy equivalence $K\approx L$ is
necessarily simple because $\text{Wh}(\pi_1(K))=0$ [18].

By an {\it open} $n$--{\it thickening} of a 
finite connected CW--complex $K$ we mean an
open smooth manifold that is diffeomorphic to
the interior of a thickening of $K$.

\subheading{Example~7.1 (Totally geodesic)}
Let $N$ be a complete manifold of nonpositive curvature
of dimension $>4$ 
that is homotopy equivalent to a finite
cell complex $K$ of dimension $\le\dim(N)-3$.
Assume that $N$ contains a (possibly noncompact) 
totally geodesic embedded submanifold $M$ of dimension $\le\dim(N)-3$ 
such that the inclusion $M\to N$ is a homotopy equivalence.

Then $N$ is an open thickening of $K$.
(Indeed, the exponential map identifies $N$ with the total space
of the normal bundle of $M$ in $N$.
Then, according to [35], the manifold $N$ has exactly
one end that is $\pi_1$-stable and
the natural homomorphism of
the fundamental group at infinity into $\pi_1(N)$
is an isomorphism.
By [25] $N$ is simple homotopy equivalent to 
a simplicial subcomplex $K^\prime$ of dimension $\le\dim(K)\le\dim(N)-3$. 
Hence, $N$ is diffeomorphic to the regular 
neighborhood of $K^\prime$ in $N$ [35].
Thus, $N$ is an open thickening of $K$.)

\subheading{7.3.~Existence of thickenings}
For each vector bundle $\xi$ over $K$ of rank $n>\max\{5,\dim(K)\}$
there exists an $n$-thickening $L$ of $K$ and a simple
homotopy equivalence
$f:K\to L$ such that $f^*TL\cong\xi$ [38, 5.1].

\subheading{7.4.~Uniqueness of thickenings}
Any tangential homotopy equivalence of $n$-thickenings
is homotopic to a diffeomorphism provided $n>5$.
[38, 5.1].

\subheading{7.5.~Uniqueness of open thickenings}
Any tangential homotopy equivalence of open $n$-thickenings
is homotopic to a diffeomorphism provided $n>4$. 
[29, pp.226--228]

\proclaim{7.6.~Theorem} 
Let $K$ be a finite, connected, 
aspherical CW-complex and 
let $n$ be an integer with $n>\max\{4,2\dim(K)\}$.
Suppose that the class $\Cal M_{a, b, \pi_1(K), n}$ 
breaks into finitely
many tangential homotopy types.

Then, for any $a\le b<0$,
the set of diffeomorphism classes of open thickenings of $K$
that belong to $\Cal M_{a, b, \pi_1(K), n}$ is \bf finite.
\endproclaim
\demo{Proof} Apply $7.5$.
\enddemo

\proclaim{7.6.~Theorem} 
Let $K$ be a finite, connected, 
aspherical CW-complex and 
let $n$ be an integer with $n>\max\{5,2\dim(K)\}$.
Suppose that the class $\Cal M_{a, b, \pi_1(K), n}$ 
breaks into finitely
many tangential homotopy types.

Then, for any $a\le b<0$,
the set of diffeomorphism classes of thickenings of $K$
whose interiors belong to $\Cal M_{a, b, \pi_1(K), n}$ is \bf finite.
\endproclaim
\demo{Proof} Apply $7.4$.
\enddemo

\proclaim{7.7.~Corollary} 
Let $K$ be a finite, connected, 
aspherical CW-complex and 
let $n$ be an integer with $n>\max\{4,2\dim(K)\}$.
Suppose that either

$\bullet$ any nilpotent subgroup of $\pi_1(K)$ has cohomological
dimension $\le 2$, or

$\bullet$ $\pi$ is not virtually nilpotent and $\pi$
does not split as a nontrivial amalgamated product or
an HNN-extension 
over a virtually nilpotent group.

Then, for any $a\le b<0$,
the set of diffeomorphism classes of open thickenings of $K$
that belong to $\Cal M_{a, b, \pi_1(K), n}$ is \bf finite.
\endproclaim
\demo{Proof} Combine $5.7$, $6.7$, and $7.6$.
\enddemo

\proclaim{7.8.~Theorem} 
Let $K$ be a finite, connected, 
aspherical CW-complex and 
let $n$ be an integer with $n>\max\{4,2\dim(K)\}$.
Suppose $\rho_k$ is a sequence of free isometric 
actions of $\pi_1(K)$ on Hadamard $n$-manifolds $X_k$
such that, for some $p_k\in X_k$,
$(X_k, p_k, \rho_k)$ is precompact in both
pointwise convergence topology and equivariant pointed
Lipschitz topology.
Assume that for each $k$, the manifold $X_k/\rho_k(\pi_1(K))$
is an open thickening of $K$.

Then the set $\{X_k/\rho_k(\pi_1(K))\}$ breaks into finitely many
diffeomorphism types.
\endproclaim

\demo{Proof} Since $K$ is a finite complex, the group
$\widetilde{KO}(K)$ is finitely generated. 
Hence $5.5$ implies that the set of manifolds $\{X_k/\rho_k(\pi_1(K))\}$ 
falls into finitely many tangential homotopy types.
According to $7.6$, 
$\{X_k/\rho_k(\pi_1(K))\}$ breaks into finitely many
diffeomorphism types.
\qed\enddemo

\head
\S 8. Finiteness for convex-cocompact groups 
\endhead

First, we recall some basic facts on convex-cocompact groups
that are well known in the constant negative curvature case
(see [11] for more details). 

Let $X$ be a Hadamard manifold with sectional curvatures
pinched between two negative constants and let $\Gamma$ be
a discrete subgroup of the isometry group of $X$.
Then $\Gamma$ acts by homeomorphisms on the ideal boundary  
$\partial_\infty X$ of $X$.
The set of points $\Omega(\Gamma)\subset\partial_\infty X$
where $\Gamma$ acts properly discontinuously
is called the {\it domain of discontinuity} of $\Gamma$.
Its complement
$\Lambda(\Gamma)=\partial_\infty X\setminus\Omega(\Gamma)$
is called the {\it limit set} of $\Gamma$. 

Fix any $\epsilon >0$.
Let $C_\epsilon(\Gamma)$ be the closed $\epsilon$-neighborhood
of the convex hull of $\Lambda(\Gamma)$ in $X$. 

There is a $\Gamma$-equivariant homeomorphism of $X\cup\Omega(\Gamma)$
and $C_\epsilon$ defined as follows. 
The fibers of orthogonal projection 
$p:X\to C_\epsilon(\Gamma)$ are 
geodesic rays orthogonal to 
$\partial C_\epsilon(\Gamma)$. 
Any point of $\Omega(\Gamma)$ is the endpoint of
such a ray and no two of these rays have the same endpoint.
Thus, the map $p$ extends to an equivariant continuous map 
$\bar{p}:X\cup\partial_\infty X\to 
C_\epsilon(\Gamma)\cup\Lambda(\Gamma)$
which is the identity on the limit set.
Contracting along the rays defines an equivariant
homeomorphism of $X\cup\partial_\infty X$ and 
$C_{2\epsilon}(\Gamma)\cup\Lambda(\Gamma)$.
which descends to a homeomorphism of 
$X\cup\Omega(\Gamma))/\Gamma$ and 
$C_{2\epsilon}(\Gamma)/\Gamma$. 

We say that $\Gamma$ is {\it convex-cocompact}
if the quotient $X\cup\Omega(\Gamma)/\Gamma$ is compact. 

Two convex--cocompact groups $\Gamma_1$ and $\Gamma_2$ 
are called {\it topologically equivalent}
if there exists a homeomorphism 
$h:X_1\cup\partial_\infty X_1\to X_2\cup\partial_\infty X_2$
that is equivariant with respect to a certain isomorphism 
of $\Gamma_1$ and $\Gamma_2$.

Given negative reals $a\le b$, a torsion-free group $\pi$,
and an integer $n$, define a class of convex-cocompact groups 
$\Cal{CC}_{a,b,\pi, n, \pi_1(\Omega)=1}$ as follows.

A convex-cocompact group $\Gamma$ of isometries
of a Hadamard manifold $X$ is said to belong to
$\Cal{CC}_{a,b,\pi, n, \pi_1(\Omega)=1}$ provided
the following three conditions hold

$\bullet$ $\Gamma$ is isomorphic to $\pi$;

$\bullet$ $\dim(X)=n$ and the sectional curvature of $X$ is within $[a,b]$;

$\bullet$ the domain of discontinuity $\Omega(\Gamma)$
of $\Gamma$ is simply-connected.

\proclaim{Proposition.~8.1}
Let $\Gamma$ be a torsion--free convex--cocompact subgroup
of the isometry group of a Hadamard manifold $X$ of dimension $>5$
such that $\Omega(\Gamma)$ is simply-connected.
Assume that $\Gamma$ is the fundamental group of a finite 
aspherical CW--complex $K$ of dimension $\le\dim(X)-3$.

Then $(X\cup\Omega(\Gamma))/\Gamma$ is a thickening of $K$.
\endproclaim
\demo{Proof}
Since $\Gamma$ be a torsion--free convex--cocompact group, 
the quotient $(X\cup\Omega(\Gamma))/\Gamma$ is a compact
manifold with boundary that is homotopy equivalent to $K$.

The homotopy equivalence is necessarily simple because
$\text{Wh}(\Gamma)=0$. (Farrell and Jones [18]
proved the vanishing of the Whitehead group of the fundamental group of
any complete manifold of pinched negative curvature.)

Finally, since $\Omega(\Gamma)$ is simply-connected, the inclusion
$\Omega(\Gamma)/\Gamma\to (X\cup\Omega(\Gamma))/\Gamma$ induces
a $\pi_1$--isomorphism.
\qed\enddemo

\proclaim{Proposition.~8.2}
Let $K$ be a finite, connected, 
aspherical CW-complex and 
let $n$ be an integer with $n>2\dim(K)$ and $n>5$.
Given $i\in\{1,2\}$, 
let $\Gamma_i$ be a torsion--free subgroup
of the isometry group of a Hadamard manifold $X_i$
such that $\Gamma_i\in\Cal{CC}_{a,b,\pi_1(K), n, \pi_1(\Omega)=1}$.

If $X_1/\Gamma_1$ and $X_2/\Gamma_2$
are tangentially homotopy equivalent,
then $\Gamma_1$ and $\Gamma_2$ are
topologically equivalent. 
\endproclaim

\demo{Proof}
Since $C_\epsilon(\Gamma_i)/\Gamma_i$
is homeomorphic to $X\cup\Omega(\Gamma_i)$,
the proposition~$8.1$ implies that 
$C_\epsilon(\Gamma_i)/\Gamma_i$ is a 
thickening of $K$.

Moreover, $C_\epsilon(\Gamma_i)/\Gamma_i$
is a codimension zero submanifold of $X/\Gamma_i$,
hence the homotopy equivalence 
$C_\epsilon(\Gamma_i)/\Gamma_i\hookrightarrow X/\Gamma_i$
is tangential.
Thus, $C_\epsilon(\Gamma_1)/\Gamma_1$ and
$C_\epsilon(\Gamma_2)/\Gamma_2$
are tangentially homotopy equivalent thickenings,
and hence they are diffeomorphic.

The diffeomorphism of compact manifolds
$C_\epsilon(\Gamma_1)/\Gamma_1\to
C_\epsilon(\Gamma_2)/\Gamma_2$ 
is necessarily bilipschitz.
Hence, it lifts to an equivariant
bilipschitz diffeomorphism $d:C_\epsilon(\Gamma_1)\to
C_\epsilon(\Gamma_2)$.

Note that $C_\epsilon(\Gamma_i)$ 
is a Gromov hyperbolic space with ideal boundary
$\Lambda(\Gamma_i)$.
Therefore, $d$ extends to an 
equivariant homeomorphism 
$C_\epsilon(\Gamma_1)\cup\Lambda(\Gamma_1)\to 
C_\epsilon(\Gamma_2)\cup\Lambda(\Gamma_2)$ [15, p35].
Finally, using equivariant homeomorphisms of
$X\cup\partial_\infty X$ and 
$C_{\epsilon}(\Gamma_i)\cup\Lambda(\Gamma_i)$,
we produce a topological equivalence of 
$\Gamma_1$ and $\Gamma_2$.
\qed\enddemo

\proclaim{8.3.~Theorem} 
Let $K$ be a finite, connected, 
aspherical CW-complex and 
let $n$ be an integer with $n>2\dim(K)$ and $n>5$.

Then the class  $\Cal{CC}_{a,b,\pi_1(K), n, \pi_1(\Omega)=1}$
falls into \bf finitely many \it
topological equivalence classes. 
\endproclaim

\demo{Proof} 
Any convex--cocompact group $\Gamma$ is word-hyperbolic
because it acts isometrically and cocompactly on a negatively curved 
space $C_\epsilon(\Gamma)$.
Hence the result follows from $6.8$ and $8.2$.
\qed 
\enddemo

\head
\S 9. Locally symmetric nonpositively curved 
manifolds up to tangential homotopy equivalence.
\endhead

For completeness we present a proof of the following result.

\proclaim{Theorem~9.1}
Let $\pi$ be a finitely presented torsion--free group
and let $X$ be a nonpositively curved symmetric space.

Then the class of manifolds of the form $X/\rho(\pi)$,
where $\rho\in\text{Hom}(\pi, \text{Isom}(X))$ is a faithful 
discrete representation, falls into finitely many
tangential homotopy types.
\endproclaim
\demo{Proof} We can assume that $\text{Hom}(\pi, \text{Isom}(X))$ 
contains a faithful discrete representation.
Let $K$ be the corresponding quotient manifold.

First note that, for any two faithful discrete
representations $\rho_1$ and $\rho_2$ that lie in the same
connected component of the analytic variety
$\text{Hom}(\pi_1(K),\text{Isom}(X))$,
we have $\tau(\rho_1)\cong\tau(\rho_2)$.
Indeed, by the covering homotopy
theorem the $X$-bundles $\tilde{K}\times_{\rho_1}X$ and
$\tilde{K}\times_{\rho_2}X$ over $K$ are isomorphic.
In particular, the pullbacks to $K$ of the vertical bundles
$\tilde{K}\times_{\rho_1}TX$ and $\tilde{K}\times_{\rho_2}TX$ 
are isomorphic as desired.

Thus, if the analytic variety $\text{Hom}(\pi,\text{Isom}(X))$
has finitely many connected components, the set 
of the manifolds of the form $X/\rho(\pi)$ where
$\rho$ is discrete and faithful
falls into finitely many tangential equivalence classes.
This is the case if
$\pi$ is finitely presented and $X$ is a nonpositively
curved symmetric space.
Indeed, represent $X$ as a Riemannian product $Y\times\Bbb R^k$
where $Y$ is a nonpositively curved symmetric space
without Euclidean factors. 
By the de Rham's theorem this decomposition is unique,
so $\text{Isom}(X)\cong\text{Isom}(Y)\times\text{Isom}(\Bbb R^k)$.
The group $\text{Isom}(Y)$ is semisimple with trivial center,
hence the analytic variety $\text{Hom}(\pi,\text{Isom}(Y))$
has finitely many connected components [21, p.567].
The same is true for $\text{Hom}(\pi,\text{Isom}(\Bbb R^k))$
because $\text{Isom}(\Bbb R^k)$ is real algebraic [21, p.567].
Hence the analytic variety
$$\text{Hom}(\pi,\text{Isom}(X))\cong
\text{Hom}(\pi,\text{Isom}(Y))\times\text{Hom}(\pi,\text{Isom}(\Bbb R^k))$$
has finitely many connected components.
\qed\enddemo

\proclaim{9.2.~Theorem} 
Let $K$ be a finite, connected, 
aspherical CW-complex and 
let $n$ be an integer with $n>\max\{4,2\dim(K)\}$.

Then the set of diffeomorphism classes of open $n$-thickenings of $K$
that admit complete locally symmetric metrics of
nonpositive sectional curvature is \bf finite.
\endproclaim
\demo{Proof} 
Note that there exist only finitely many symmetric Hadamard manifolds
of a given dimension. Hence, the result follows from
$7.5$ and $9.1$.
\qed\enddemo

\proclaim{9.3.~Theorem} 
Let $K$ be a finite, connected, 
aspherical CW-complex and let $X$ be
a symmetric negatively curved Hadamard $n$-manifold
with $n>\max\{5,2\dim(K)\}$.

Let $\rho_1$ and $\rho_2$ be injective representations
of $\pi_1(K)$ into the isometry group of $X$
that lie in the same connected component of the
representation variety 
$\text{Hom}(\pi_1(K),\text{Isom}(X))$.
Suppose that, for $i\in\{1,2\}$,
$\rho_i(\pi_1(K))$ is a convex-cocompact group
with simply-connected domain of discontinuity.

Then $\rho_1$ and $\rho_2$ are conjugate be a 
homeomorphism of $X\cup\partial_\infty X$
\endproclaim
\demo{Proof}
Recall that the are four kinds of symmetric spaces of
negative sectional curvature, namely, they are
hyperbolic spaces over the reals, complex numbers,
quaternions and Cayley numbers.
The spaces have sectional curvatures
pinched between $-4$ and $-1$.
Applying $8.1$ we conclude that 
$X\cup\Omega(\rho_i(\pi_1(K)))/\rho_i(\pi_1(K))$
is a thickening of $K$ for $i=1,2$.

It follows from the proof of $9.1$ that
the homotopy equivalence induced by
$\rho_1\circ(\rho_2)^{-1}$
is tangential because $\rho_1$ and $\rho_2$ 
lie in the same connected component.
Thus, by $7.4$, the homotopy equivalence
is homotopic to a diffeomorphism that lifts
to a smooth conjugacy of $\rho_1$ and $\rho_2$
on $X\cup\Omega(\rho_1(\pi_1(K)))$ and 
$X\cup\Omega(\rho_2(\pi_1(K)))$. 
Repeating the argument of $8.2$, we deduce
that the conjugacy extends to an equivariant
self-homeomorphism of $X\cup\partial_\infty X$. 
\qed\enddemo

\enddocument